\documentclass[12pt]{article}
\usepackage{hyperref}
\usepackage{amsfonts}
\usepackage{enumerate}
\usepackage{amsthm}
\usepackage{amsmath}
\usepackage{graphicx}

\begin{document}

\begin{center}
{\large\bf A mathematical model for plasticity and damage: A discrete calculus formulation}

\vskip.20in

Ioannis Dassios$^{1,3,4}$ Andrey Jivkov$^{1}$ Andrew Abu-Muharib$^{1,2}$ Peter James$^{2}$\\[2mm]
{\footnotesize
$^{1}$Mechanics of Physics of Solids Research Team, School of Mechanical Aerospace \& Civil Engineering,
The University of Manchester, UK\\[5pt]
$^{2}$AMEC Foster Wheeler, Birchwood, Risely, Warrington, UK\\[5pt]
$^{3}$MACSI, Department of Mathematics \& Statistics, University of Limerick, Ireland\\[5pt]
$^{4}$ERC, Electricity Research Centre, University College Dublin, Ireland}
\end{center}

{\footnotesize
\noindent
\textbf{Abstract.} In this article we propose a discrete lattice model to simulate the elastic, plastic and failure behaviour of isotropic materials. Focus is given on the mathematical derivation of the lattice elements, nodes and edges, in the presence of plastic deformations and damage, i.e. stiffness degradation. By using discrete calculus and introducing non-local potential for plasticity, a force-based approach, we provide a  matrix formulation necessary for software implementation. The output 
is a non-linear system with allowance for elasticity, plasticity and damage in lattices.  This is the key tool for explicit analysis of micro-crack generation and population growth in plastically deforming metals, leading to macroscopic degradation of their mechanical properties and fitness for service. An illustrative example, analysing a local region of a node, is given to demonstrate the system performance. \\\\[3pt]
{\bf Keywords} : discrete calculus, lattice model, steel microstructure, plasticity, damage, non-linear system
\\[3pt]

\vskip.2in
\section{Introduction}
Lattice models for analysis of deformation and fracture of solids have been developed over the last thirty years mainly for quasi-brittle materials, such as concretes and rocks \cite{Cus}, \cite{Gri}, where local behaviour is elastic-brittle, i.e. the only mechanism for energy dissipation is the generation of new surfaces (brittle micro-cracking). Lattice models contain set of sites connected by bonds; in the language of algebraic topology this is a 1-complex embedded in $\mathbb{R}^2$ or $\mathbb{R}^3$, i.e. $2D$ or $3D$ graphs. The challenge with lattices is that they intend to represent a \textit{continuous} solid by a discrete system. Specifically, the stored energy in any lattice region is required to be equivalent to the stored energy in the corresponding continuum region. This is used to derive a link between properties of lattice elements, e.g. bond stiffness coefficients, and macroscopic properties of the material, e.g. elastic modulus and Poisson's ratio. For irregular graphs this is not possible exactly; approximate calibrations of lattice properties for known macroscopic properties are used, e.g. \cite{Yip}. For graphs with some regularity, the link can be established rigorously. For example, isotropic materials can be represented exactly by $2D$ graphs based on hexagonal structure, e.g. \cite{Kar}, and by $3D$ graphs based on truncated octahedral structure \cite{Zha}. Other possible regular $3D$ graphs can represent cubic elasticity, i.e. three parameter elastic behaviour, but not isotropic elasticity, i.e. two parameter behaviour, typical for most engineering materials \cite{Wa}. Nevertheless, lattices are being used because of the need to represent failure in materials, albeit not always with exactly calibrated local properties. Failure is a generation of a new internal surface, i.e. a non-topological concept. Hence, the classical continuum mechanics does not work, as it is a thermodynamic \textit{bulk} theory.

Our work addresses cleavage fracture, which is a phenomenon of fast and catastrophic crack propagation. It can be observed in all metals with body-centred cubic crystal lattices due to the smaller number of slip planes, i.e. reduced ability to dissipate energy by plastic deformation. We are specifically interested in ferritic steels, used for example to manufacture nuclear reactor pressure vessels, where cleavage is a potential cause for serious concern. Cleavage is known to be triggered by rupturing of brittle second-phase particles, such as carbides, which typically decorate the grain boundaries of the polycrystal and can rupture due to plastic overload from the surrounding grains. The existing modelling of this phenomenon, the so called local approach to cleavage fracture, is based on the weakest-link statistic, i.e. global failure probability derives from the failure probabilities of individual particles, following prescribed size distribution but treated as independent events \cite{Ber}. This approach works well for predicting cleavage fracture at very low temperatures, where plasticity is limited, the number of micro-cracks formed prior to cleavage is small and they are spatially scattered, i.e. only a set of largest particles failed. The probability of cleavage fracture, however, decreases rapidly with the increase of temperature, where plasticity is enhanced. This apparent increase of material toughness with temperature cannot be captured with the existing modelling strategy, despite of the many improvements in the particle rapture criterion over the years \cite{Jam}, \cite{Pin}. However, accurate assessment with reduced conservatism is needed for more economic exploitation of reactors, including better planning of inspection intervals and life-extension decisions. From the current status of modelling it can be deduced that the plasticity enhancement with temperature leads to generation of increasingly large micro-crack populations, in which case the argument behind the weakest-link statistic is violated, e.g. there is a significant effect of micro-crack interactions prior to ultimate cleavage. This will explain the reduced probability of cleavage, hence increased toughness, of the material with temperature. Lattice models are particularly suitable to study the generation and interaction of micro-cracks as they grow, coalesce and progress to final failure. To this end we present a lattice development relevant to plastically deforming metals. 

A mathematically rigorous treatment of lattices can be achieved when they are analysed as graphs \cite{Gra}. However, the discrete exterior calculus (DEC) specialized to graphs in this reference, is developed for scalar problems. This means that the nodal unknown (a 0-cochain) is scalar, i.e. temperature, pressure, concentration, etc. The gradient of this is also a scalar field over the edges (1-cochain). Things go relatively simply for such physical problems. In mechanics this approach is rather more difficult as we have a vector-valued nodal field, namely a displacement vector assigned to each node. Article \cite{Ya} is the only one to our knowledge which attempts to apply DEC to mechanical problems in elastic settings. Our aim is to build upon this theoretical basis and present a graph-theoretical approach to lattices including elastic, plastic and damage behaviour.
\\\\
\textit{Microstructure representation}
\\\\
Micro-structures of metals and alloys develop by a crystallization process, which starts at spatially randomly distributed nuclei and finishes when neighbouring growing crystals touch each other to form grain boundaries \cite{And1}. Under uniform temperature distribution, the individual crystals grow with equal rates and therefore the final micro-structure is described by the Voronoi diagram constructed around the set of nuclei \cite{And2}. The $3D$ Voronoi diagram is constructed by the intersection of planes bisecting normally the segments connecting each pair of neighbouring nuclei, as illustrated in Figure 1. Thus, each crystal is a polyhedron representing the neighbourhood of a nucleus, where all points are closer to the nucleus than to any other nucleus. 
\begin{figure}[h]
\begin{center}
\includegraphics[width=0.80\textwidth]{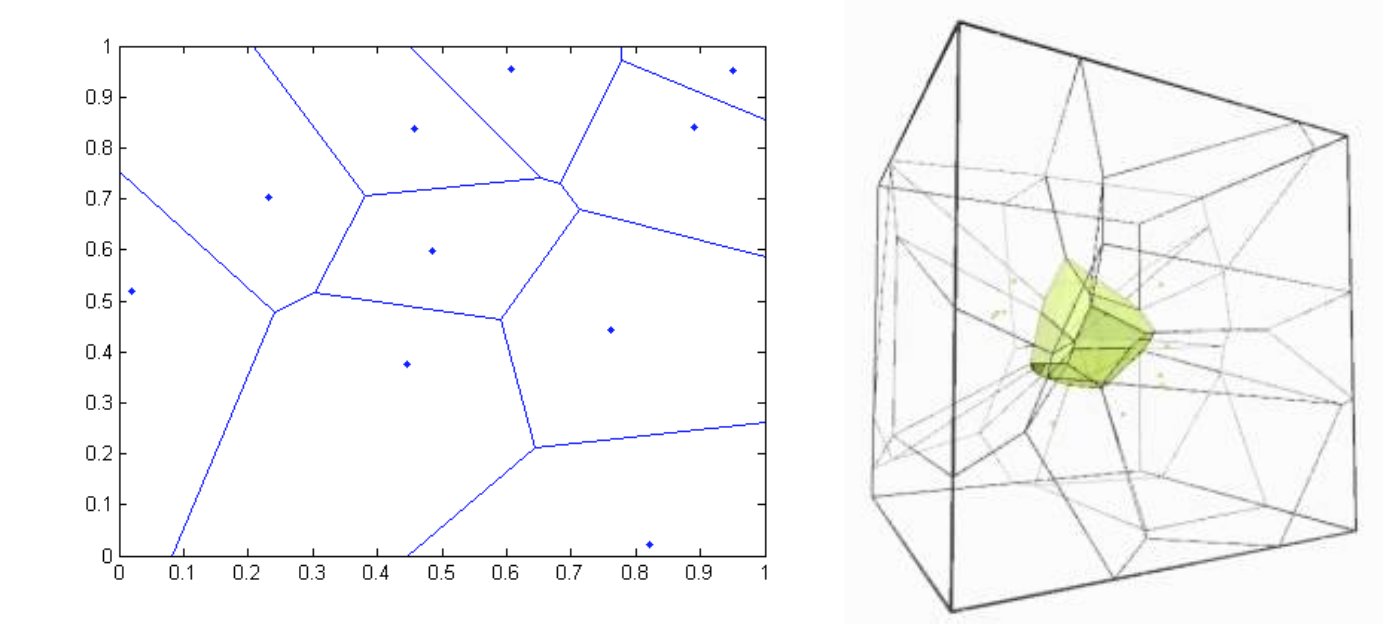}
\caption{Illustration of Voronoi diagrams constructed around a set of points: 2D construction (left) providing comprehensive picture of the process; 3D construction (right) with the polyhedron formed around the central point shown in green.}
\label{figure1}
\end{center}
\end{figure}
Initially, we seek to represent the poly-crystal with a regular tessellation of the $3D$ space, shown in Figure 2, where a grain is topologically equivalent to the average grain in arbitrary Voronoi tessellation \cite{Jiv}. This simplifies the geometry and allows for exact calibration of the lattice emerging from the tessellation \cite{Zha}.
\begin{figure}[h]
\begin{center}
\includegraphics[width=0.80\textwidth]{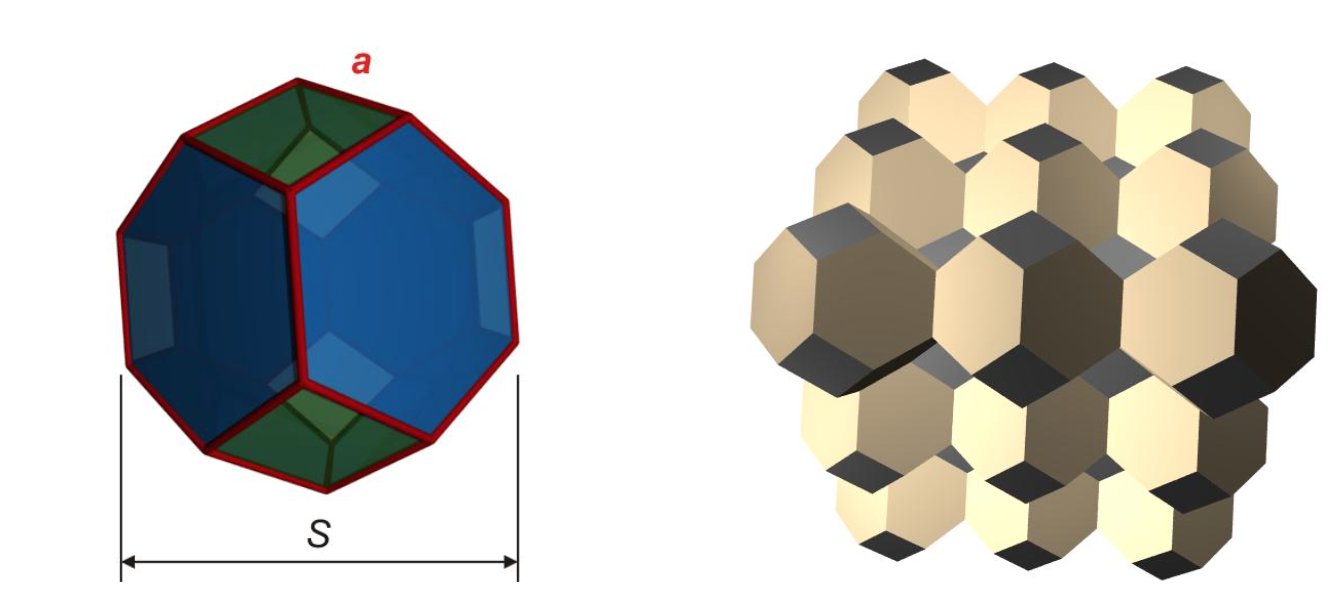}
\caption{Truncated octahedron (left), representing an average single grain. Assembly of cells, filling compactly a 3D region (right), representing topologically averaged polycrystal.}
\label{figure2}
\end{center}
\end{figure}
\\\\
\textit{Cell complexes and their duals}\\\\
Irrespective of whether we use an irregular tessellation of space (Voronoi around nuclei) or the regular tessellation of Figure 2, the material is subdivided into compactly packed and non-intrusive $3D$ cells. Using terminology from algebraic topology, see \cite{Ha}, the polycrystal is a $p$-complex, where $p = 3$ embedded in $\mathbb{R}^d$, where $d = 3$. Thus, each crystal is a 3-cell, bounded by a number of 2-cells (planar polygonal faces), each of which is bounded by a number of 1-cells (edges), each of which is bounded by two 0-cells (nodes).
It should be noted that $p = d$ is not a requirement for $p$-complexes embedded in $\mathbb{R}^d$; $p$-complexes with $1\leq p \leq d$ can be formed and treated similarly. For example, the collection of all faces of the regularized polycrystal (or grain boundaries of the real microstructure) without the cell (grain) interiors is a 2-complex embedded in $\mathbb{R}^3$. It contains 0-cells, the vertices of the truncated octahedrons (or quadruple points in the real microstructure), 1-cells, the edges of the truncated octahedrons (or triple lines), and 2-cells, the faces of the truncated octahedrons (or grain boundaries). The structure of thin-walled closed cell foam can be viewed as a 2-complex. Similarly, the collection of all edges of the regularized polycrystal (or triple lines in the real microstructure) without the cell (grain) interiors and faces (boundaries) is a 1-complex embedded in $\mathbb{R}^3$. It contains 0-cells, the vertices of the truncated octahedrons (or quadruple points in the real microstructure), and 1-cells, the edges of the truncated octahedrons (or triple lines). The structure of thin-ligament open cell foam can be viewed as a 1-complex. Illustration of the three possibilities for complexes embedded in $\mathbb{R}^3$ is given in Figure 3. The historical name for 1-complex is a graph, in this particular case the graph of edges in the right figure.
\begin{figure}[h]
\begin{center}
\includegraphics[width=1.00\textwidth]{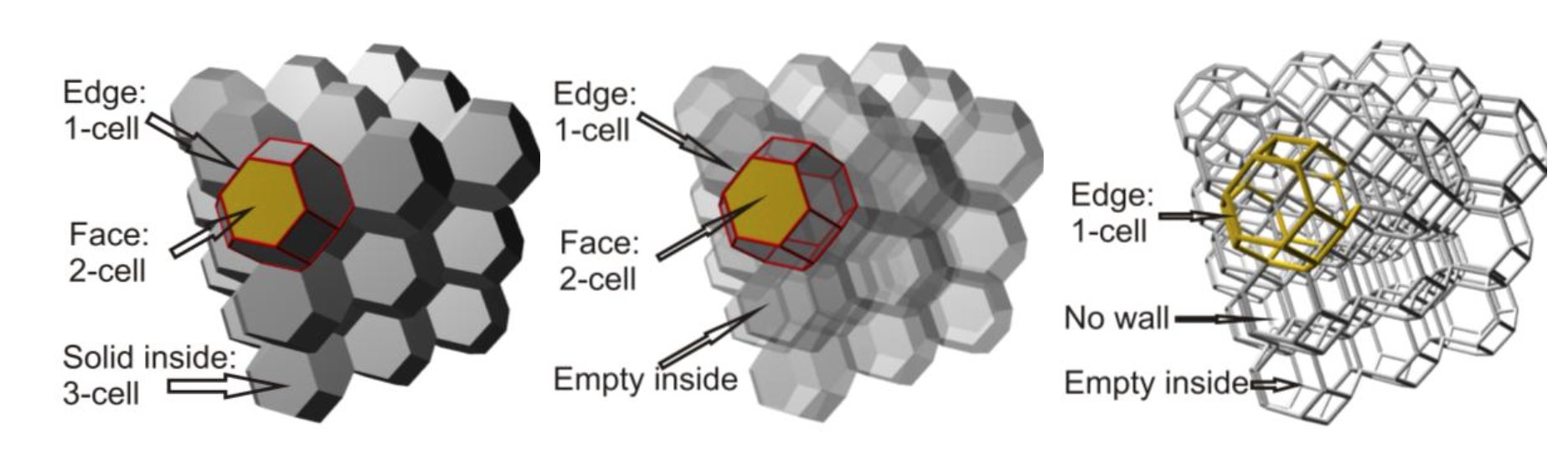}
\caption{Illustration of 3-complex (left), 2-complex (centre), and 1-complex (right), embedded in $\mathbb{R}^3$. Illustration based on the topologically averaged polycrystal.}
\label{figure3}
\end{center}
\end{figure}

Important for the analysis on complexes is that each $p$-complex can be used to build a dual $p$-complex, where the $s$-cells of the dual correspond to $(p-s)$-cells of the original or primal, for each $0 \leq s \leq p$. Specifically for a 3-complex, the 0-cells of the dual correspond to 1-cells of the primal, the 1-cells of the dual correspond to 2-cells of the primal, the 2-cells of the dual correspond to 1-cells of the primal, and the 3-cells of the dual correspond to 0-cells of the primal. For a 3-complex embedded in $3D$, such as the material microstructure, this means that the nodes of the dual will correspond to the grains, the edges of the dual will correspond to grain boundaries, the faces of the dual will correspond to triple lines, and the cells of the dual will correspond to quadruple points.

Considering the representation of the material as a 3-complex, the topology of the site-bond model is derived as a 1-complex (graph) by placing sites (graph nodes) at the centres of all 3-cells and bonds (graph edges) between neighbouring sites as illustrated in Figure 4. Each site (node) has 14 incident bonds (edges): eight normal to hexagonal faces (nearer neighbours) denoted by $B_1$ with length $L_1 = \frac{\sqrt{3}}{2}S=\sqrt{6}a$ and six normal to square faces (further neighbours) denoted by $B_2$ with length $L_2 = S = \sqrt{8}a$. We keep the expressions with both $S$ and $a$, because $S$ is more convenient when the 3-complex is sized using an average grain volume, $S = (2V_{GR})^\frac{1}{3}$, while $a$ is more convenient when it is sized using an average triple line length, $a = L_{TL}$, in the real microstructure. Note, that the construction of the graph follows the construction of a 3-complex dual to the one representing the microstructure, but it terminates before introducing 2-cells dual to the 1-cells (edges) and 3-cells dual to the 0-cells (vertices) of the original 3-complex.
\begin{figure}[h]
\begin{center}
\includegraphics[width=1.00\textwidth]{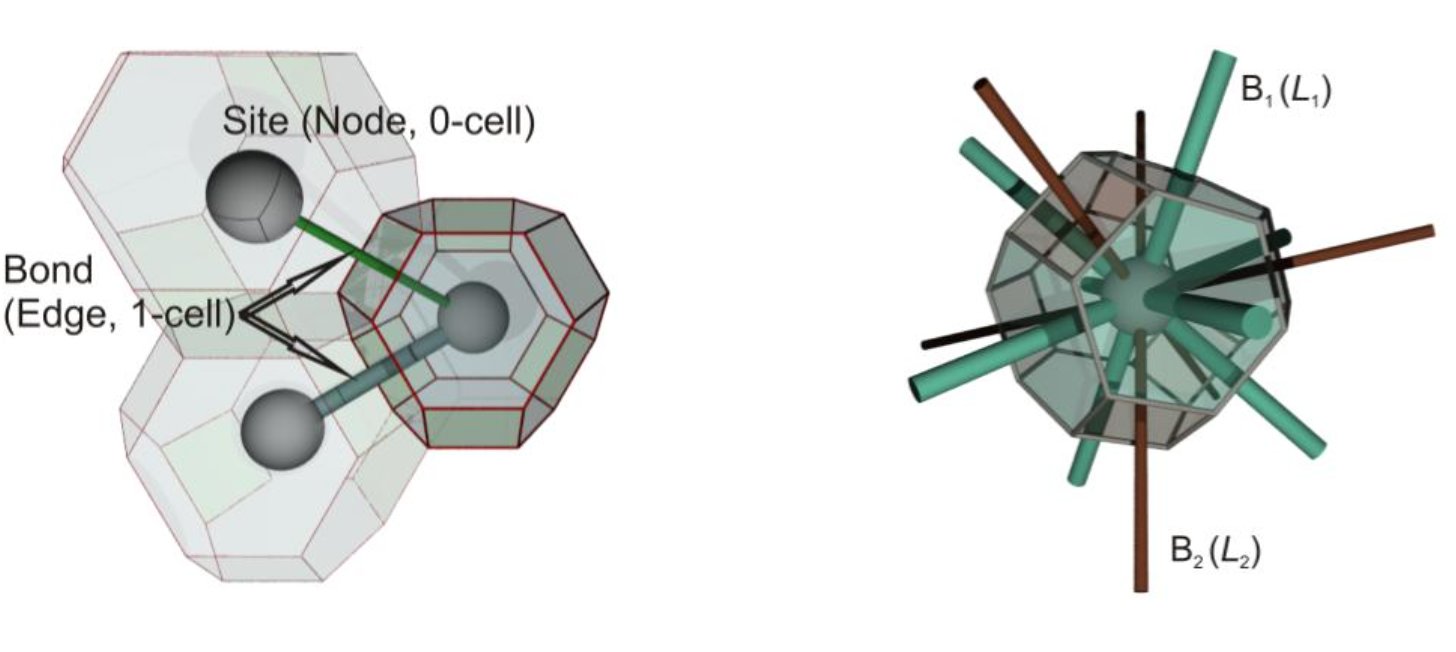}
\caption{Illustration of site-bond model construction as a reduction of 3-complex to 1- complex (graph). Left figure shows a 3-cell with two neighboring 3-cells which goes into 3 sites (0-cells, graph nodes) connected by bonds (1-cells, graph edges). Right figure shows a single 3-cell with the 14 bonds (1-cells, graph edges) in the substituting site-bond topology. The two types of bonds, denoted by $B_1$ and $B_2$, have distinct lengths $L_1$ and $L_2$.}
\label{figure4}
\end{center}
\end{figure}

\section{The model with calculus on discrete manifolds}

In the discrete calculus theory the boundary of a $p$-complex is represented by the so-called incidence matrix, see \cite{Gra}. Specifically for 1-complexes (graphs), the incidence matrix, takes the simplest form denoted by $A=[a_{ij}]^{j=1,2,...,n}_{i=1,2,...,m}$. Where 
\[
a_{ij}=\left\{\begin{array}{cc} 0,&\textnormal{if node $j$ is not a node of edge $i$}\\
1,&\textnormal{if node $j$ is the first node of edge $i$}\\-1,&\textnormal{if node $j$ is the second node of edge $i$}\end{array}\right\}.
\]
For the the node-edge incidence matrix $A\in\mathbb{R}^{m\times n}$, where $m$ is the number of edges and $n$ is the number of nodes of the graph. Row $i$, corresponding to edge $i$, contains exactly one 1 in the column corresponding to the first node of the edge and one -1 in the column corresponding to the second node of the edge. All remaining coefficients of the matrix are zeros. The incidence matrix contains not only the topology of the graph, but also represent the boundary operator.  Thus application of $A\in\mathbb{R}^{m\times n}$ to a nodal function is an edge function - the gradient of the nodal. Further, application of the transposed incidence matrix to an edge function is a nodal function - the divergence of the edge. Since we will know the initial positions of the nodes and the edges from a partiocular cellular structure, the incidence matrix will be known. 

Assume a lattice of $n$ nodes connected through $m$ edges. Let $D_i\in\mathbb{R}^3$, $i=1,2,...,n$ be a discrete (vector-valued) function over nodes and $b_i\in\mathbb{R}^3$, $i=1,2,...,m$ a discrete function over edges. The relation
\[
AD=b.
\]
states that $b$ is the gradient of $D$, where 
\[
D=\left[\begin{array}{c}D_1\\D_2\\\vdots\\D_n \end{array}\right]\in\mathbb{R}^{n\times 3}, b=\left[\begin{array}{c}b_1\\b_2\\\vdots\\b_m \end{array}\right]\in\mathbb{R}^{m\times 3}.
\]

Let 
\[
X=\left[\begin{array}{c}X_1\\X_2\\\vdots\\X_n \end{array}\right]\in\mathbb{R}^{n\times 3} ,y=\left[\begin{array}{c}y_1\\y_2\\\vdots\\y_m \end{array}\right]\in\mathbb{R}^{m\times 3},
\]
where $X_i\in\mathbb{R}^3$, $i=1,2,...,n$ are the nodal coordinate and $y_i\in\mathbb{R}^3$, $i=1,2,...,m$ the edge coordinates, the modulus of which provide edge lengths. Then
\begin{equation}\label{eq1}
AX=y.
\end{equation}
In \eqref{eq1}, the components of $y$ are unknown, but some components of $X$ are known from the essential boundary conditions required to complete the boundary value problem. The loading on our system will be applied only via essential boundary conditions, i.e. displacement control, where the coordinates of some nodes change from initial to prescribed final values.

Solid materials accommodate strain from external loading by reversible (elastic) rearrangement, giving rise to internal stresses, and by dissipating energy via slip (plasticity) or separation (surface generation). In discrete settings it has been shown that the stress between two regions is a traction vector always normal to their interface \cite{Ya}. In the graph framework, this means that the stress is a vector (force, not a tensor as in continuum mechanics) acting along the current (deformed) orientation of an edge.
Let $n_i=\frac{y_i}{\left|y_i\right|}$, $i=1,2,...,m$,  be the unit vectors along edges, where $\left|y_i\right|$, $i=1,2,...,m$, are the edge lengths after some deformation. The edge forces $F_i$, $i=1,2,...,m$, are therefore given by
\[
F_i=\left|F_i\right|n_i=\frac{\left|F_i\right|}{\left|y_i\right|}y_i, \quad \forall i=1,2,...,m.
\]
This can be summarised by 
\begin{equation}\label{eq2}
F=K(y)y,
\end{equation}
where
\begin{equation}\label{eq2a}
F=\left[\begin{array}{c}\left|F_1\right|\\\left|F_2\right|\\\vdots\\\left|F_m\right| \end{array}\right]\in\mathbb{R}^m,\quad K(y)=diag\left\{\frac{\left|F_1\right|}{\left|y_1\right|},\frac{\left|F_2\right|}{\left|y_2\right|},\dots.\frac{\left|F_m\right|}{\left|y_m\right|}\right\}\in\mathbb{R}^{m\times m}.
\end{equation}

Let $\left|b_i\right|$, $i=1,2,...,m$, be the initial lengths of bonds. Then $\left|F_i\right|$, $i=1,2,...,m$, are related to the edge elongations, $\left|b_i\right|-\left|y_i\right|$, via potentially a non-smooth function as illustrated in Figure 5.
\begin{figure}[h]\label{1}
\begin{center}
\includegraphics[width=0.70\textwidth]{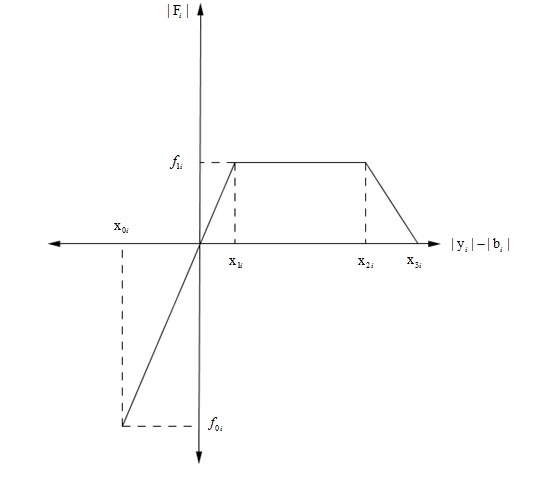}
\caption{Conceptual edge forces-displacement relation.}
\label{1}
\end{center}
\end{figure}

Bond elongation, positive x-axis in Figure 5, results in positive or tensile force (or zero after $x_{3i}$). Bond contraction, results in negative or compressive force. Bond elongation is related to the established in mechanics notions of stretch and strain. Stretch is the ratio between deformed and initial lengths, i.e. $\frac{\left|y_i\right|}{\left|b_i\right|}$, and strain is defined as the natural logarithm of stretch in the general case. For stretches around 1, the strain is approximated by the engineering strain $\frac{\left|y_i\right|-\left|b_i\right|}{\left|b_i\right|}$. Mathematically the stretch can be anything from zero to one under contraction and from one to infinity under elongation. In practice, the possible stretches depend on the material in question. Rubber for example allows for stretches more than 2 (more than 100 percent elongation) without any plasticity or failure. Steels become plastic at strains between 0.2 percent and 0.5 percent. Hence the difference $\frac{\left|y_i\right|-\left|b_i\right|}{\left|b_i\right|}$ is really small at the onset of plasticity, represented by the line with zero slope in Figure 5. Surface separation, represented by the line with negative slope (softening branch), is initiated at several percent strain and ultimate failure occurs typically between 40 percent and 60 percent.

Let $B_i$, $i=1,2,...,n$ be external forces at nodes, given as natural boundary conditions. These can all be zero as in the case of displacement-controlled loading, but can in general represent force-controlled loading. The balance of linear momentum at nodes (force equilibrium) requires that 
\begin{equation}\label{eq3}
A^TF=B.
\end{equation}
where $A^T \in\mathbb{R}^{n\times m}$ is the transpose of the incidence matrix $A$ and 
\[
B=\left[\begin{array}{c}B_1\\B_2\\\vdots\\B_n \end{array}\right]\in\mathbb{R}^{n\times 3}.
\]
In summary, at each node we have a position vector (3 components in three coordinate directions) and a force vector (3 components in the same directions). The discrete position field on nodes is $X$, containing $3n$ numbers, where $n$ is the number of nodes. The discrete force field on nodes $B$, contains $3n$ numbers. Applying $AX$ gives us the discrete position field on edges $y$. This is a field containing $3m$ numbers, three coordinate differences per edge. The position field on edges gives rise to a force field in the edges, which resists the deformation. Thus $F$ is a discrete force field on the edges with $3m$ components. Hence we have $F=f(y)=K(y)y$, which is a material-dependent function defining the resistance to deformation. Applying $A^TF$ gives us a discrete force field on the nodes. This must equal $B$ to balance the linear momentum in the system, i.e. ensure equilibrium. 

We note the exclusion principle for nodal fields in the system definition: an element of $B$ is known (the force in a given direction at a given node) if and only if the corresponding element of $X$ is unknown (the displacement in the same direction at the same node). This means that when a nodal force is prescribed (a natural boundary condition), the system solution provides its conjugate displacement. Inversely, when a nodal displacement is prescribed (essential boundary condition), the system solution provides its conjugate force.

Substitution of \eqref{eq2} into \eqref{eq3} gives
\[
A^TK(y)y=B,
\]
whereby using \eqref{eq1} into the above expression provides the formal description of the mechanical deformation of the graph with
\begin{equation}\label{eq4}
A^TK(AX)AX=B.
\end{equation}

\section{Solving the non-linear system}

System \eqref{eq4} is a non-linear system. This is because of the diagonal matrix $K$ defined in \eqref{eq2a}. We provide the following Theorem
\\\\
\textbf{Theorem 3.1.} Consider the system \eqref{eq4}. Then 
\begin{enumerate}[(a)]
\item An effective linearization of system \eqref{eq4} is
\begin{equation}\label{eq5}
\tilde{A}X=B.
\end{equation} 
Where $\tilde{A}=A^T\tilde{K}A$ and for $0<\epsilon<<1$
\[
\tilde{K}=
\left\{\begin{array}{cc}
diag[\frac{f_{0i}}{2x_{0i}}(1-\frac{\left|b_i\right|}{\left|b_i\right|+x_{0i}})\pm\epsilon]_{1\leq i\leq m},
&x_{0i}\leq\left|y_i\right|-\left|b_i\right|\leq 0\\
diag[\frac{f_{1i}}{2x_{1i}}(1-\frac{\left|b_i\right|}{\left|b_i\right|+x_{1i}})\pm \epsilon]_{1\leq i\leq m},
&0\leq\left|y_i\right|-\left|b_i\right|\leq x_{1i}\\
diag[\frac{f_{1i}(x_{2i}-x_{1i})}{2(\left|b_i\right|+x_{1i})(\left|b_i\right|+x_{2i})}\pm \epsilon]_{1\leq i\leq m},
&x_{1i}\leq\left|y_i\right|-\left|b_i\right|\leq x_{2i}\\
diag[\frac{f_{1i}(x_{3i}-\left|b_i\right|)}{2(x_{3i}+\left|b_i\right|)(x_{2i}+\left|b_i\right|)}\pm \epsilon]_{1\leq i\leq m},
&x_{2i}\leq\left|y_i\right|-\left|b_i\right|\leq x_{3i}
\end{array}\right\}.
\]

\item Let 
\[
\tilde{A}=
\left[
\begin{array}{cc}\tilde{A}_{11}&\tilde{A}_{12}\\
\tilde{A}_{21}&\tilde{A}_{22}\end{array}\right].
\]
Then system \eqref{eq5} can be divided into the following subsystems
\begin{equation}\label{eq6}
\tilde{A}_{11}X_p=B_p-\tilde{A}_{12}X_q
\end{equation}
and
 \begin{equation}\label{eq7}
B_q=\tilde{A}_{13}X_p+\tilde{A}_{22}X_q.
\end{equation} 
From the above systems only \eqref{eq6} has to be solved. Then $X_p$ can be replaced in \eqref{eq7} and $B_q$ is easily computed.
\item For the solution of \eqref{eq6}
\begin{enumerate}[(i)]
\item If $\tilde{A}_{12}$ is full rank, then
\begin{equation}\label{eq8}
X_p=\tilde{A}_{11}^{-1}(B_p-\tilde{A}_{12}X_q).
\end{equation}
\item If $\tilde{A}_{11}$ is rank deficient, then an optimal solution $\hat{X}_p$ for \eqref{eq6} is given by
\begin{equation}\label{eq9}
   \hat X_{p}=(\tilde{A}_{11}^T\tilde{A}_{11}+E^TE)^{-1}\tilde{A}_{11}^T(B_p-\tilde{A}_{12}X_q).
\end{equation}
Where $E$ is a matrix such that $\tilde{A}_{11}^T\tilde{A}_{11}+E^TE$ is invertible and $\left\|E\right\|_2=\theta$, $0<\theta<<1$. Where $\left\|\cdot\right\|_2$ is the Euclidean norm and with $()^T$ we denote the transpose tensor.
\end{enumerate}
\end{enumerate}
\textbf{Proof.} For the proof of (a) we consider the system \eqref{eq4}, \eqref{eq2a}. $\forall i=1,2,...,m$ we will seek bounds for $\frac{\left|F_i\right|}{\left|y_i\right|}$. We have the following cases. \\\\
For $x_{0i}\leq\left|y_i\right|-\left|b_i\right|\leq 0$ we have $\left|F_i\right|=\frac{f_{0i}}{x_{0i}}(\left|y_i\right|-\left|b_i\right|)$, see Figure 5. Or, equivalently,
\begin{equation}\label{eq10}
\frac{\left|F_i\right|}{\left|y_i\right|}=\frac{f_{0i}}{x_{0i}}(1-\frac{\left|b_i\right|}{\left|y_i\right|}).
\end{equation}
Since $x_{0i}\leq\left|y_i\right|-\left|b_i\right|\leq 0$, we have $x_{0i}+\left|b_i\right|\leq\left|y_i\right|\leq \left|b_i\right|$, or, equivalently, $\frac{x_{0i}+\left|b_i\right|}{\left|b_i\right|}\leq\frac{\left|y_i\right|}{\left|b_i\right|}\leq 1$, or, equivalently, $1\leq\frac{\left|b_i\right|}{\left|y_i\right|}\leq \frac{\left|b_i\right|}{x_{0i}+\left|b_i\right|}$ or, equivalently,
\[
1-\frac{\left|b_i\right|}{x_{0i}+\left|b_i\right|}\leq1-\frac{\left|b_i\right|}{\left|y_i\right|}\leq 0. 
\]
Hence, by using the above inequality and \eqref{eq10} we get
\[
\frac{\left|F_i\right|}{\left|y_i\right|}\cong\frac{f_{0i}}{2x_{0i}}(1-\frac{\left|b_i\right|}{\left|b_i\right|+x_{0i}})\pm\epsilon.
\]
For $0\leq\left|y_i\right|-\left|b_i\right|\leq x_{1i}$ we have $\left|F_i\right|=\frac{f_{1i}}{x_{1i}}(\left|y_i\right|-\left|b_i\right|)$, see Figure 5, or, equivalently,
\begin{equation}\label{eq11}
\frac{\left|F_i\right|}{\left|y_i\right|}=\frac{f_{1i}}{x_{1i}}(1-\frac{\left|b_i\right|}{\left|y_i\right|}).
\end{equation}
Since $0\leq\left|y_i\right|-\left|b_i\right|\leq x_{1i}$, we have 
\[
0\leq1-\frac{\left|b_i\right|}{\left|y_i\right|}\leq 1-\frac{\left|b_i\right|}{x_{1i}+\left|b_i\right|}. 
\]
Hence, by combining the above inequality with \eqref{eq11} we get
\[
\frac{\left|F_i\right|}{\left|y_i\right|}\cong\frac{f_{1i}}{2x_{1i}}(1-\frac{\left|b_i\right|}{\left|b_i\right|+x_{1i}})\pm \epsilon.
\]
For $x_{1i}\leq\left|y_i\right|-\left|b_i\right|\leq x_{2i}$ we have $\left|F_i\right|=f_{1i}$, see Figure 1, or equivalently,
\begin{equation}\label{eq12}
\frac{\left|F_i\right|}{\left|y_i\right|}=\frac{f_{1i}}{\left|y_i\right|}.
\end{equation}
Since $x_{1i}\leq\left|y_i\right|-\left|b_i\right|\leq x_{2i}$, we have $x_{1i}+\left|b_i\right|\leq\left|y_i\right|\leq x_{2i}+\left|b_i\right|$, or equivalently
\[
\frac{1}{x_{2i}+\left|b_i\right|}\leq\frac{1}{\left|y_i\right|}\leq \frac{1}{x_{1i}+\left|b_i\right|}
\]
or, equivalently, by using the above inequality and \eqref{eq12} we get
\[
\frac{f_{1i}}{x_{2i}+\left|b_i\right|}\leq\frac{\left|F_i\right|}{\left|y_i\right|}\leq \frac{f_{1i}}{x_{1i}+\left|b_i\right|}. 
\]
Hence
\[
\frac{\left|F_i\right|}{\left|y_i\right|}\cong\frac{f_{1i}(x_{2i}-x_{1i})}{2(\left|b_i\right|+x_{1i})(\left|b_i\right|+x_{2i})}\pm \epsilon.
\]
For $x_{2i}\leq\left|y_i\right|-\left|b_i\right|\leq x_{3i}$ we have $\left|F_i\right|=\frac{f_{1i}}{x_{2i}-x_{3i}}(\left|y_i\right|-\left|b_i\right|)+\frac{x_{3i}f_{1i}}{x_{3i}-x_{2i}}$, see Figure 5, or equivalently,
\begin{equation}\label{eq13}
\frac{\left|F_i\right|}{\left|y_i\right|}=\frac{f_{1i}}{x_{2i}-x_{3i}}(1-\frac{\left|b_i\right|}{\left|y_i\right|})+\frac{x_{3i}f_{1i}}{x_{3i}-x_{2i}}\frac{1}{\left|y_i\right|}.
\end{equation}
Since $x_{2i}\leq\left|y_i\right|-\left|b_i\right|\leq x_{3i}$, we have $x_{2i}+\left|b_i\right|\leq\left|y_i\right|\leq x_{3i}+\left|b_i\right|$, or equivalently, $\frac{x_{2i}+\left|b_i\right|}{\left|b_i\right|}\leq\frac{\left|y_i\right|}{\left|b_i\right|}\leq \frac{x_{3i}+\left|b_i\right|}{\left|b_i\right|}$, or equivalently 
\[
\frac{\left|b_i\right|}{x_{3i}+\left|b_i\right|}\leq\frac{\left|b_i\right|}{\left|y_i\right|}\leq \frac{\left|b_i\right|}{x_{2i}+\left|b_i\right|},
\]
or, equivalently,
\[
1-\frac{\left|b_i\right|}{x_{2i}+\left|b_i\right|}\leq1-\frac{\left|b_i\right|}{\left|y_i\right|}\leq 1-\frac{\left|b_i\right|}{x_{3i}+\left|b_i\right|}. 
\]
Similarly,
\[
\frac{1}{x_{3i}+\left|b_i\right|}\leq\frac{1}{\left|y_i\right|}\leq \frac{1}{x_{2i}+\left|b_i\right|}.
\]
Hence, by using the above two inequalities combined with \eqref{eq13} we have 
\[
\frac{\left|F_i\right|}{\left|y_i\right|}\cong \frac{f_{1i}}{2(x_{2i}-x_{3i})}(\frac{\left|b_i\right|}{\left|b_i\right|+x_{2i}}-\frac{\left|b_i\right|}{\left|b_i\right|+x_{3i}})+\frac{x_{3i}f_{1i}}{2(x_{3i}-x_{2i})}(\frac{1}{\left|b_i\right|+x_{2i}}-\frac{1}{\left|b_i\right|+x_{3i}})\pm\epsilon,
 \]
 or, equivalently,
\[
\frac{\left|F_i\right|}{\left|y_i\right|}\cong \frac{f_{1i}\left|b_i\right|(x_{3i}-1)}{2(\left|b_i\right|+x_{2i})(\left|b_i\right|+x_{2i})}\pm\epsilon.
\]
For the proof of (b) system \eqref{eq5} can be written as
\[
\left[\begin{array}{cc}\tilde{A}_{11}&\tilde{A}_{12}\\\tilde{A}_{21}&\tilde{A}_{22}\end{array}\right]
\left[\begin{array}{c}X_p\\X_q\end{array}\right]=
\left[\begin{array}{c}B_p\\B_q\end{array}\right],
\]
or, equivalently,
\[
\begin{array}{c}\tilde{A}_{11}X_p+\tilde{A}_{12}X_q=B_p,\\\tilde{A}_{13}X_p+\tilde{A}_{22}X_q=B_q.\end{array}
\]
From the above expressions we get the subsystems \eqref{eq6}, \eqref{eq7}.\\\\
For the proof of (c), for (i), since $\tilde{A}_{11}$ is full rank, we arrive easy at \eqref{eq8}. For (ii), since $\tilde{A}_{11}$ is rank deficient, if $[B_p-\tilde{A}_{12}X_q]\notin colspan \tilde A_{11}$ system \eqref{eq6} has no solutions and if $[B_p-\tilde{A}_{12}X_q]\in colspan \tilde A_{11}$ system \eqref{eq6} has infinite solutions. Let 
\[
\hat L(\hat X_p)=\hat L+E\hat X_p,
\]
such that the linear system 
\[
\tilde A_{11}\hat X_p=\hat L(\hat X_p),
\]
or, equivalently the system
\[
(\tilde A_{11}-E)\hat X_p=\hat L
\]
has a unique solution. Where $E$ is a matrix such that $\tilde A_{11}^T\tilde A_{11}+E^TE$ is invertible, $\left\|E\right\|_2=\theta$, $0<\theta<<1$ and $E\hat X_p$ is orthogonal to $\hat L -\tilde A_{11}\hat X_p$.  Hence we want to solve the following optimization problem
\[
\begin{array}{c}min\left\|(B_p-\tilde{A}_{12}X_q)-\hat L\right\|_2^2\\s.t.\quad (\tilde A_{11}-E)\hat X_p=\hat L,\end{array}
\]
or, equivalently,
\[
min\left\|(B_p-\tilde{A}_{12}X_q)-(\tilde A_{11}-E)\hat X_p\right\|_2^2.
\]
or, equivalently,
\[
min\left\{\left\|(B_p-\tilde{A}_{12}X_q)-\tilde A_{11}\hat X_p\right\|_2^2+\left\|E\hat X_p\right\|_2^2\right\}.
\]
The reason for using the matrix $E$ is because the matrix $\tilde A_{11}$ is rank deficient and hence the matrix $\tilde A_{11}^T\tilde A_{11}$ is singular and not invertible. To sum up, we seek a solution $\hat X_p$ minimizing the functional
\[
H_1(\hat X_p)=\left\|B_p-\tilde{A}_{12}X_q-\tilde A_{11}\hat X_p\right\|_2^2+\left\|E\hat X_p\right\|_2^2.
\]
Expanding $H_1(\hat X_p)$ gives
\[
H_1(\hat X_p)=(B_p-\tilde{A}_{12}X_q-\tilde A_{11}\hat X_p)^T(B_p-\tilde{A}_{12}X_q-\tilde A_{11}\hat X_p)+(E\hat X_p)^TE\hat X_p,
\]
or, equivalently,
\[
H_1(\hat X_p)=(B_p-\tilde{A}_{12}X_q)^T(B_p-\tilde{A}_{12}X_q)-2(B_p-\tilde{A}_{12}X_q)^T\tilde A_{11}\hat X_p+(\hat X_p)^T\tilde A_{11}^T\tilde A_{11}\hat X_p+(\hat X_p)^TE^TE\hat X_p
\]
because $(B_p-\tilde{A}_{12}X_q)^T\tilde A_{11}\hat X_p=(\hat X_p)^T\tilde A_{11}^T(B_p-\tilde{A}_{12}X_q)$. Furthermore
\[
\frac{\partial}{\partial \hat X_p}H_1(\hat X_p)=-2\tilde A_{11}^T(B_p-\tilde{A}_{12}X_q)+2\tilde A_{11}^T\tilde A_{11}\hat X_p+2E^TE\hat X_p.
\]
Setting the derivative to zero, $\frac{\partial}{\partial \hat X_p}H_1(\hat X_p)=0$, we get
\[
(\tilde A_{11}^T\tilde A_{11}+E^TE)\hat X_p=\tilde A_{11}^T(B_p-\tilde{A}_{12}X_q).
\]
The solution is then given by
\[
\hat X_p=(\tilde A_{11}^T\tilde A_{11}+E^TE)^{-1}\tilde A_{11}^T(B_p-\tilde{A}_{12}X_q).
\]
Hence an optimal solution of \eqref{eq6} is given by \eqref{eq9}. The proof is completed.

\section{Numerical example}

\begin{figure}[h]\label{1}
\begin{center}
\includegraphics[width=0.50\textwidth]{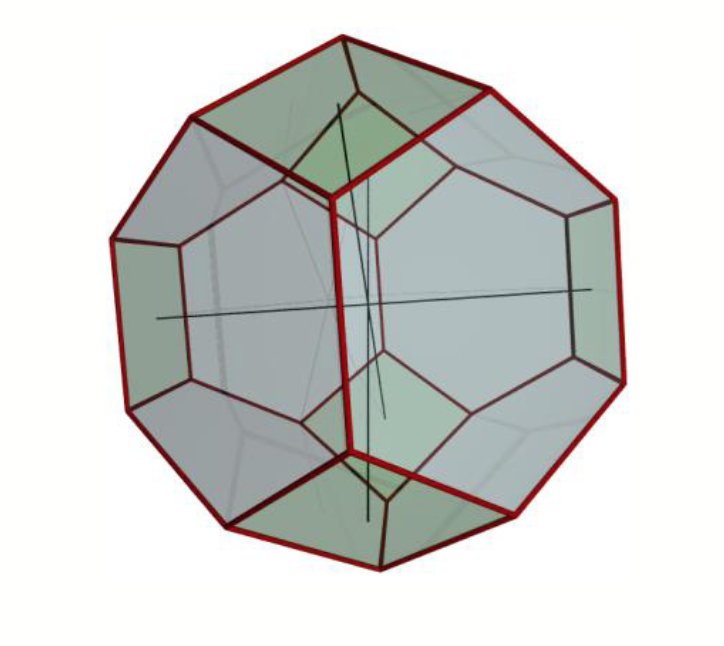}
\caption{A unit cell}
\label{figure6}
\end{center}
\end{figure}

We assume a unit cell of 15 nodes and 14 edges. Let $D_i$, $i=1,2,...,15$, $b_i$ and $\left|b_i\right|$, $x_{0i}$, $x_{1i}$, $x_{2i}$, $x_{3i}$, $f_{0i}$, $f_{1i}$, $i=1,2,...,14$, for $i=1,2,...,14$, be given by Table 1.
\begin{table}[ht] 
\centering 
\begin{tabular}{c cc ccc cccc} 
\hline\hline 
$i$&$D_i$&$b_i$&$\left|b_i\right|$&$x_{0i}$&$x_{1i}$&$x_{2i}$&$x_{3i}$&$f_{0i}$&$f_{1i}$ \\ [0.5ex] 
\hline 
1&(0,0,0)&(0,-1,1)&$\sqrt 2$& -0.1&0.1&0.3&0.4&-0.1&0.1\\ 
2&(0,1,0)&(0,1,0)&  1&              -0.2&0.2&0.3&0.5&-0.1&0.1\\
3&(0,0,1)&(1,-1,-1)&  $ \sqrt 3$&             -0.1&0.1&0.2&0.3&-0.2&0.2\\
4&(0,1,1)&(0,1,0)& 1&               -0.1&0.1&0.3&0.4&-0.1&0.1\\
5 &(1,0,0)&(0,-1,1)&   $\sqrt 2$ &           -0.1&0.1&0.3&0.4&-0.1&0.1\\ 
6&(1,1,0)&(0,1,0)&     1  &         -0.2&0.2&0.3&0.5&-0.1&0.1\\
7&(1,0,1)&(1,-1,-1)&    $\sqrt 3$&           -0.1&0.1&0.2&0.3&-0.2&0.2\\
8&(1,1,1)&(0,1,0)&    1           & -0.1&0.1&0.3&0.4&-0.1&0.1\\
9&(2,0,0)&(0,-1,1)&     $\sqrt 2$  &         -0.2&0.2&0.3&0.5&-0.1&0.1\\
10&(2,1,0)&(0,1,0&       1  &      -0.1&0.1&0.2&0.3&-0.2&0.2\\
11&(2,0,1)&(0,1,0)&       1  &       -0.1&0.1&0.3&0.4&-0.1&0.1\\
12&(2,1,1)&(0,-1,1)&       $\sqrt 2$ &        -0.2&0.2&0.3&0.4&-0.2&0.2\\
13&(2,2,1)&(0,1,0)&        1  &      -0.1&0.1&0.3&0.4&-0.1&0.1\\ 
14&(2,1,2)&(-2,-1,-2) &       3&        -0.1&0.1&0.3&0.4&-0.1&0.1\\
15&(2,2,2)&&&&&&&&
\\[1ex] 
\hline 
\end{tabular} 
\caption {Nodes, edges and other given values}
\end{table}
Then
\[
A=\left[\begin{array}{cccccccc}0&0 & 0 &0&0 &\dots&0  & 0  \\0&-1  & 1 &0&0 &\dots&0  & 0  \\0&0 & -1  & 1& 0&  \dots&0  & 0  \\\vdots  & \vdots  & \vdots  &\vdots   &  \vdots  &  \ddots&  \vdots&\vdots   \\0& 0  & 0 &0&0 &\ldots& -1& 1 \end{array}\right].
\]

In addition, if $X=\left[\begin{array}{c}X_p\\X_q\end{array}\right]\in\mathbb{R}^{15\times 3}$, $B=\left[\begin{array}{c}B_p\\B_q\end{array}\right]\in\mathbb{R}^{15\times 3}$, with $p=7$, $q=8$, let $X_p$, $B_q$ be unknown and $X_q$, $B_p$ given as follows:
\begin{enumerate}[(a)]
\item For $x_{0i}\leq\left|y_i\right|-\left|b_i\right|\leq 0$
\[
X_q=\left[\begin{array}{ccc}0.1&1.1&0.1\\0.1&0.1&1.1\\ 0.1&1.1&1.1\\ 1.1&0.1&0.1\\1.1&1.1&0.1\\1.1&0.1&1.1\\1.1&1.1&1.1\\1.1&0.1&0.1 \end{array}\right],
B_p=\left[\begin{array}{ccc}0.1&1.1&0.1\\0.1&0.1&1.1\\ 0.1&1.1&1.1\\ 1.1&0.1&0.1\\1.1&1.1&0.1\\1.1&0.1&1.1\\1.1&1.1&1.1 \end{array}\right];
\]
\item For $0\leq\left|y_i\right|-\left|b_i\right|\leq x_{1i}$
\[
X_q=\left[\begin{array}{ccc}0.2&1.2&0.2\\0.2&0.2&1.2\\ 0.2&1.2&1.2\\ 1.2&0.2&0.2\\1.2&1.2&0.2\\1.2&0.2&1.2\\1.2&1.2&1.2\\1.2&0.2&0.2 \end{array}\right],
B_p=\left[\begin{array}{ccc}0.2&1.2&0.2\\0.2&0.2&1.2\\ 0.2&1.2&1.2\\ 1.2&0.2&0.2\\1.2&1.2&0.2\\1.2&0.2&1.2\\1.2&1.2&1.2 \end{array}\right];
\]
\item For $x_{1i}\leq\left|y_i\right|-\left|b_i\right|\leq x_{2i}$
\[
X_q=\left[\begin{array}{ccc}0.3&1.3&0.3\\0.3&0.3&1.3\\ 0.3&1.3&1.3\\ 1.3&0.3&0.3\\1.3&1.3&0.3\\1.3&0.3&1.3\\1.3&1.3&1.3\\1.3&0.3&0.3 \end{array}\right],
B_p=\left[\begin{array}{ccc}0.3&1.3&0.3\\0.3&0.3&1.3\\ 0.3&1.3&1.3\\ 1.3&0.3&0.3\\1.3&1.3&0.3\\1.3&0.3&1.3\\1.3&1.3&1.3 \end{array}\right];
\]
\item For  $x_{2i}\leq\left|y_i\right|-\left|b_i\right|\leq x_{3i}$
\[
X_q=\left[\begin{array}{ccc}0.4&1.4&0.4\\0.4&0.4&1.4\\ 0.4&1.4&1.4\\ 1.4&0.4&0.4\\1.4&1.4&0.4\\1.4&0.4&1.4\\1.4&1.4&1.4\\1.4&0.4&0.4 \end{array}\right],
B_p=\left[\begin{array}{ccc}0.4&1.4&0.4\\0.4&0.4&1.4\\ 0.4&1.4&1.4\\ 1.4&0.4&0.4\\1.4&1.4&0.4\\1.4&0.4&1.4\\1.4&1.4&1.4 \end{array}\right].
\]
\end{enumerate}
From Theorem 3.1 the solution of systems \eqref{eq6}, \eqref{eq7}, will be
\begin{enumerate}[(a)]
\item For $x_{0i}\leq\left|y_i\right|-\left|b_i\right|\leq 0$
\[
X_p=\left[\begin{array}{ccc} 
            0    &     0     &    0\\
111.2090&  110.9590&  152.1000\\
  111.0840&  110.8340&  150.7250\\
  107.6199&   90.0494&  112.6199\\
  101.1199&   83.5494&  101.1199\\
   84.1493&   66.5788&   84.1493\\
   79.7743 &  63.4538&   79.7743
 \end{array}\right],
B_q=\left[\begin{array}{ccc}   
-4.6000 &  -3.4000&   -4.8000\\
    0.0000&   -0.7657 &   0.2000\\
   -0.1000 &   0.6657 &   0.1000\\
    0.1000  & -0.3000 &  -0.1000\\
         0&    0.4828&   -0.2828\\
         0   &-0.4828 &   0.2828\\
         0  &  0.2667 &   0.0667\\
         0   &-0.0667&   -0.0667
 \end{array}\right].
\]

\item For $0\leq\left|y_i\right|-\left|b_i\right|\leq x_{1i}$
\[
X_p=\left[\begin{array}{ccc} 
              0&         0 &        0\\
 -129.8438& -127.5938& -170.7349\\
 -129.5938& -127.3438& -169.2349\\
 -122.6656& -103.0951& -127.6656\\
 -114.6656&  -95.0951& -114.6656\\
  -94.8666&  -75.2961 & -94.8666\\
  -89.8666&  -71.5461&  -89.8666
 \end{array}\right],
B_q=\left[\begin{array}{ccc}   
 -5.2000&   -4.4000&   -5.0000\\
         0&    0.7657&   -0.2000\\
    0.1000&   -0.6657&   -0.1000\\
   -0.1000  &  0.3000&    0.1000\\
    0.0000&   -0.4828&    0.2828\\
         0 &   0.4828&   -0.2828\\
         0 &  -0.2667&   -0.0667\\
         0 &   0.0667 &   0.0667

 \end{array}\right];
\]
 
\item For $x_{1i}\leq\left|y_i\right|-\left|b_i\right|\leq x_{2i}$

\[
X_p=\left[\begin{array}{ccc} 
     0&         0&         0\\
    0.3591&   1.3559&    0.3749\\
    0.3581 &   1.3549&    0.3708\\
    0.3564 &   1.3504&    0.3634\\
    0.3431 &   1.3371&    0.3431\\
    0.3308  &  1.3248  &  0.3308\\
    0.3164  &  1.3136 &   0.3164
 \end{array}\right],
B_q=\left[\begin{array}{ccc}   

         -5.8&        138.2&       -148.8\\
            0  &    -696.42&          143\\
         -132  &     685.42&          132\\
          132&         -275    &     -132\\
   0&       419.71   &   -276.71\\
            0&      -419.71 &      276.71\\
  0&         1166      &   1023\\
            0  &      -1023 &       -1023

 \end{array}\right];
\]

\item For  $x_{2i}\leq\left|y_i\right|-\left|b_i\right|\leq x_{3i}$
\[
X_p=\left[\begin{array}{ccc} 
        0  &          0 &           0\\
    -0.026822 &       0.986 &    -0.12908\\
    -0.021694      &0.99113    & -0.11113\\
    0.0074869      & 1.0568  & -0.0089966\\
     0.043751      &  1.093     &0.043751\\
      0.10245       &1.1517      &0.10245\\
      0.16655       & 1.203      &0.16655
 \end{array}\right],
B_q=\left[\begin{array}{ccc}   

         -6.4 &     -66.067 &      54.267\\
   0&       132.45&      -60.667\\
       22.286&      -94.071     & -22.286\\
      -22.286  &     82.952       &22.286\\
   0&     -91.33      & 30.664\\
   0&        91.33     & -30.664\\
            0&      -146.97      &-86.308\\
            0 &      86.308       &86.308

 \end{array}\right].
\]

\end{enumerate}
In addition, the edge elongations are given by Tables 2, 3.

\begin{table}[ht] 
\centering 
\begin{tabular}{c cc} 
\hline\hline 
$y_i$&   $x_{0i}\leq\left|y_i\right|-\left|b_i\right|\leq 0$&    $0\leq\left|y_i\right|-\left|b_i\right|\leq x_{1i}$ \\ [0.5ex] 
\hline 
$y_1$&(0, 0, 0)&(0, 0, 0)\\ 
$y_2$&(-0.1250, -0.1250, -1.3750)&(0.25, 0.25, 1.5)\\
$y_3$&(-3.4641, -20.7846, -38.1051)&(6.9282, 24.2487, 41.5692)\\
$y_4$&(-6.5, -6.5, -11.5)&(8 ,8, 13)\\
$y_5$&(-16.9706, -16.9706, -16.9706)&(19.799, 19.799, 19.799)\\ 
$y_6$&(-4.375, -3.125, -4.375)&(5, 3.75, 5)\\
$y_7$&(-79.6743, -62.3538, -79.6743)&(90.0666, 72.7461, 90.0666)\\
$y_8$&(0, -1, 1)&(0, -1, 1)\\
$y_9$&(0, 1, 0)&(0, 1, 0)\\
$y_{10}$&(1, -1, -1)&(1, -1, -1)\\
$y_{11}$&(0, 1, 0)&(0, 1, 0)\\
$y_{12}$&(0, -1, 1)&(0, -1, 1)\\
$y_{13}$&(0, 1, 0)&(0, 1, 0)\\ 
$y_{14}$&(0, -1, -1)&(0, -1, -1)
\\[1ex] 
\hline 
\end{tabular} 
\caption {Edge elongations}
\end{table}

\begin{table}[ht] 
\centering 
\begin{tabular}{c cc} 
\hline\hline 
$y_i$&   $x_{1i}\leq\left|y_i\right|-\left|b_i\right|\leq x_{2i}$&    $x_{2i}\leq\left|y_i\right|-\left|b_i\right|\leq x_{3i}$ \\ [0.5ex] 
\hline 
$y_1$&(0, 0, 0)&(0, 0, 0)\\ 
$y_2$&(-0.00096154, -0.00096154, -0.0041667)&(0.0051282, 0.0051282, 0.017949)\\
$y_3$&(-0.0016951, -0.0045203, -0.0073454)&(0.029181, 0.065657, 0.10213)\\
$y_4$&(-0.013287, -0.013287, -0.02028)&(0.036264, 0.036264, 0.052747)\\
$y_5$&(-0.012328, -0.012328, -0.012328)&(0.058701, 0.058701, 0.058701)\\ 
$y_6$&(-0.014423, -0.011218, -0.014423)&(0.064103, 0.051282, 0.064103)\\
$y_7$&(-0.016386, -0.013561, -0.016386)&(0.23345, 0.19697, 0.23345)\\
$y_8$&(0, -1, 1)&(0, -1, 1)\\
$y_9$&(0, 1, 0)&(0, 1, 0)\\
$y_{10}$&(1, -1, -1)&(1, -1, -1)\\
$y_{11}$&(0, 1, 0)&(0, 1, 0)\\
$y_{12}$&(0, -1, 1)&(0, -1, 1)\\
$y_{13}$&(0, 1, 0)&(0, 1, 0)\\ 
$y_{14}$&(0, -1,-1)&(0,-1,-1)
\\[1ex] 
\hline 
\end{tabular} 
\caption {Edge elongations}
\end{table}

\section*{Conclusions}

In this article we focused on the mathematical derivation of the behavior of lattice elements, nodes and edges, in the presence of non-linear deformations. We presented an efficient integrated lattice formulation and demonstrated the methodology with a simple example. An extension of this work is planned to test and verify the method, by comparing  simulations of lattice elastic-plastic-damage behavior to experimentally measured material responses. \\\\[3pt]

\section*{Acknowledgments}
I. Dassios is supported by Science Foundation Ireland (award 09/SRC/E1780). A. Abu-Muharib would like to acknowledge the support from EPSRC via the Nuclear EngD Doctoral Training Centre at The University of Manchester's Dalton Nuclear Institute and the ongoing support from sponsoring company AMEC-Clean Energy Europe.

\end{document}